\documentclass[12pt,reqno]{amsart}

\usepackage{amssymb}
\usepackage[arrow,matrix]{xy}

\newcommand{\Gal}{\text{\rm Gal}}
\newcommand{\id}{\text{\rm id}}
\newcommand{\Z}{\mathbb{Z}}

\begin{document}

\title{Hilbert~90 for Biquadratic Extensions}

\author[R.~Dwilewicz]{Roman Dwilewicz}

\author[J.~Min\'{a}\v{c}]{J\'an Min\'a\v{c}}

\author[A.~Schultz]{Andrew Schultz}

\author[J.~Swallow]{John Swallow}

\maketitle

\markboth{}{}

\newtheorem{theorem}{Theorem}
\newtheorem{corollary}{Corollary}

\theoremstyle{definition}
\newtheorem*{remark*}{Remark}

\parskip=10pt plus 5pt minus 2pt
\parindent=20pt

\noindent\textbf{1. INTRODUCTION: ``COMMERCE WITH THE GREAT
UNIVERSE''.}

\begin{verse}
    \tiny{\emph{Here something stubborn comes,\\
    Dislodging the earth crumbs\\
    And making crusty rubble.\\
    It comes up bending double,\\
    And looks like a green staple.\\
    It could be seedling maple,\\
    Or artichoke, or bean.\\
    That remains to be seen.} \\
    \qquad ---Richard Wilbur, ``Seed Leaves''}
\end{verse}

We linger over the simplest theorems, the exquisite ones, succinct
and stately. Something quickens: dare we  hope that such a theorem,
a result of already uncommon brevity, might grow into a result of
uncommon grandeur? We fancy it a seedling, a tiny seed-case somehow
containing within itself its own generalization. We imagine that its
hidden truths, unknown to us now, await only fertile soil and close
observation.  In the quadratic case of Hilbert's Theorem 90 we find
just such a seedling.

Turning it over in our hands, we naturally wonder what we might one
day behold.  Friends tell us to expect what is affectionately known
as ``Hilbert 90,'' a result in the theory of cyclic extensions.  But
quadratic Hilbert 90, as we shall now call the quadratic case, does
not lead only to Hilbert 90. Those who have ventured into Galois
cohomology know that Andreas Speiser sensed something different, a
cohomological result in a noncohomological era. Emmy Noether
recognized the power of Speiser's result and used it to great
effect: the result then became known as Noether's theorem. As we
shall see, even tiny quadratic Hilbert 90 contains suggestions of
cohomology.

More recently, others have felt fresh intimations, imagining
something different still. The spectacular announcement by Vladimir
Voevodsky of the validity of the Bloch-Kato conjecture contains, in
fact, a generalization of Hilbert's Theorem 90 to Milnor
$K$-theory.\footnote{The interested reader should consult
\cite{LeSc} and \cite{H} for some initial background on Hilbert's
Theorem~90 and \cite[p.~30]{We} for its cohomological
generalization.  To observe the use of Hilbert~90-type theorems in
the partially published work of Markus Rost and Voevodsky on the
Bloch-Kato conjecture, see \cite{Vo1} and \cite{Vo2}.  For further
original sources on Hilbert~90 and its cohomological generalization
see Ernst Kummer's early discovery of a special case \cite{Kum},
followed by Speiser's result \cite{S} and Noether's application
\cite{N}.}

In this article we propose to linger just long enough to investigate
the very first shoots of the generalization from quadratic to
biquadratic field extensions.  The quadratic case is well known and
already significant.  The biquadratic case, particularly the
difference between it and the quadratic case, appears new. Observing
the shoots break the surface and bend back upon themselves, we
encounter the possibility that the unfurling leaves aim in an
uncharted direction, somewhere between the traditional Hilbert 90
result on cyclic extensions and Speiser's and Noether's
cohomological result. Exploring this new possibility, we discern
connections among multiplicative groups of fields, values of binary
quadratic forms, a bit of module theory over group rings, and even
Galois cohomology.

\noindent \textbf{2. ``ARTICHOKE, OR BEAN''?: THE QUADRATIC SEED-CASE.}
The classical form of Hilbert 90 states that the kernel of the
norm operator on a cyclic extension is as small as possible. We
begin by making this statement precise in the context of a quadratic
extension.

Let $F$ be a field of characteristic different from two. Then any
quadratic extension $L$ of $F$ takes the form $L = F (\sqrt{a})$ for
a suitable element $a$ of $F$ that is not contained in the set $F^2$
of squares of elements of $F$. The set $\{1,\sqrt{a}\}$ is a basis
for the vector space $L$ over the field $F$, so each element $\ell$
in $F (\sqrt {a})$ can be expressed uniquely as $\ell = f_1 + f_2
\sqrt{a}$ with $f_1$ and $f_2$ in $F$.

The field extension $L/F$ has a special property: if an irreducible
polynomial $g(x)$ with coefficients in $F$ has a root in $L$, then $L$
contains a complete set of distinct roots of $g$.  Put more precisely,
if for some element $\ell_1$ of $L$ and some polynomial $h(x)$ with
coefficients in $L$ we can factor $g(x)=(x-\ell_1)\cdot h(x)$, then
there exists a constant $c$ in $L$ and a collection of distinct
elements $\ell_i$ in $L$ such that $g(x)=c\cdot \prod_i (x-\ell_i)$.
Field extensions $L/F$ enjoying this property are called \emph{Galois
  extensions}.  As a set, the \emph{Galois group} $\Gal(L/F)$ of a
Galois extension $L/F$ consists of all field automorphisms of $L$
leaving $F$ pointwise fixed.  The group operation is function
composition.

The great discovery of Galois theory is that some fundamental
properties of Galois extensions $L$ of $F$ are already captured by
the corresponding Galois groups $\Gal(L/F)$.  For instance, if the
Galois group has precisely $n$ subgroups, then Galois theory tells
us that there exist precisely $n$ fields $M$ satisfying $F\subseteq
M\subseteq L$.  To emphasize the importance of certain
group-theoretic properties of the Galois group, we use the same
adjectives to modify the corresponding Galois field extension.  For
instance, a field extension $L/F$ is called \emph{cyclic} if $L/F$
is a Galois extension and the associated Galois group $\Gal(L/F)$ is
a cyclic group.  An important result in Galois theory---and one that
we will have frequent occasion to use---is the fact that the subset
of $L$ consisting of elements fixed by every element of the Galois
group $\Gal(L/F)$ is precisely $F$.

With these definitions, we return to the quadratic extension
$L=F(\sqrt{a})$ and find that $L/F$ is a cyclic extension whose
Galois group $G = \Gal (L/F)$ contains precisely two elements: the
identity automorphism, denoted by $\id$, and $\sigma$, which sends
$\ell=f_1+f_2\sqrt{a}$ to $\sigma (\ell) = f_1 - f_2 \sqrt{a}$.  We
can verify that the elements of $L$ fixed by $\Gal(L/F)$ are, in
fact, those lying in $F$ by noticing that if $\sigma(\ell)=\ell$,
then $f_2=0$, placing $\ell=f_1$ in $F$.

The norm function $N_{L/F} : L \to F$ for such a quadratic extension
is defined by the formula $N_{L/F} (\ell) = \ell \cdot \sigma (\ell)
= f_1^2 - a f_2^2$, and, omitting zero from the domain, we find that
$N_{L/F}$ restricts to a homomorphism from the multiplicative group
$L^\times := L \setminus \{ 0 \}$ to the multiplicative group
$F^\times$ of $F$.  Since there is no chance for confusion, we write
$N_{L/F}$ for this restricted function.

It is natural to determine the kernel of this homomorphism. Since $N
(\ell) = \ell \cdot \sigma (\ell) = N (\sigma(\ell))$ for each $\ell$
in $L^\times$, we see that elements of the form $\ell/{\sigma
(\ell)}$ certainly lie in $\ker N_{L/F}$. Hilbert's Theorem~90 says
that $\ker N_{L/F}$ contains no other elements.

\begin{theorem}[\textbf{Hilbert 90 for
$\mathbf{F(\mathbf{\root\of{\mathbf{a}}})/F}$}]\label{th:quadratic H90}
    For a quadratic extension $F(\sqrt{a})$ of $F$ it is the case
    that
    \begin{equation*}
        \ker \big(N_{F(\sqrt{a})/F} : F (\sqrt{a})^
        \times \to F^\times\big) = \left\{ \frac{\ell}{\sigma
        (\ell)} : \ell \in F (\sqrt{a})^\times \right\}.
    \end{equation*}
\end{theorem}

Not surprisingly, this simple statement is already important. For
example, the classical parameterization of Pythagorean triples is a
beautiful consequence of this statement. The idea of the proof can
be traced back to Olga Taussky \cite[pp.~808--809]{T} (who was in
fact a coeditor of Hilbert's collected works), and Noam Elkies
independently discovered the proof in a short, attractive note
\cite{El}.  See also Takashi Ono's book \cite[pp.~4--5]{O}.

\begin{proof}[Proof of Theorem~\ref{th:quadratic H90}]
    Let $t$ belong to $\ker N_{L/F}$. If $t = -1$, then $-1 = t =
    {\sqrt{a}}/{\sigma (\sqrt{a})}$. Assume therefore that $t \ne
    -1$. Set $\ell = 1 + t$. Then
    \begin{equation*}
        \sigma (\ell) = 1 + \sigma (t) = \sigma (t) \cdot t + \sigma
        (t) = \sigma (t) \cdot (t+1) = \sigma (t) \cdot \ell.
    \end{equation*}
    Hence ${\sigma (\ell)}/{\ell} = \sigma (t)$. Applying
    $\sigma$ to both sides and remembering that $\sigma^2 = \id$ we
    obtain ${\ell}/{\sigma (\ell)} = t$, as required. Therefore
    \begin{equation*}
        \ker N_{L/F} \subseteq \left\{ \frac{\ell}{\sigma (\ell)} \
        :\ \ell \in L^\times \right\}.
    \end{equation*}
    We have already observed the reverse inclusion.
\end{proof}

The generalization of this result to an arbitrary cyclic extension is
well known as ``Hilbert 90.''  Many of us have come through precisely
such underbrush while ambling across our own fields.  Speiser saw,
however, that in the right environment the seedling that is quadratic
Hilbert 90 produces something similar but different---something that
our naturalist friends might call ``homologous.''

Let us review the example in Theorem~\ref{th:quadratic H90} and
observe how the seed might develop into an elegant result in Galois
cohomology. Consider a map $f:G\to L^\times$ satisfying
\begin{equation*}
    f (\gamma_1 \cdot \gamma_2) = (\gamma_1 (f
    (\gamma_2))) \cdot f(\gamma_1) \quad\quad
    (\gamma_1, \gamma_2 \in G).
\end{equation*}
Because of its similarity to a homomorphism, such a map is called a
\emph{crossed homomorphism}.  From the equality $f(\id) =
(\id(f(1))) \cdot f(\id) = f(\id)^2$ we deduce that $f(\id) = 1$.
As $G$ is especially simple in our context, it follows that a
crossed homomorphism is determined by $f(\sigma)$.  From $1 =
f(\sigma^2) = (\sigma( f(\sigma)))\cdot f(\sigma)$ we see that
$f(\sigma)$ lies in $\ker N_{L/F}$. Conversely, choosing any element
$t$ from $\ker N_{L/F}$, we can define a crossed homomorphism $f : G
\to L^\times$ by $f(1) = 1$, $f(\sigma) = t$. Appealing to
Theorem~\ref{th:quadratic H90}, we have $f(\sigma) = \ell /
\sigma(\ell)$ for some $\ell$ in $L^\times$.  In particular we find
that for any crossed homomorphism $f$ there exists $\ell$ in
$L^\times$ such that $f(g) = {\ell}/{g(\ell)}$ for all $g$ in $G$.

For an arbitrary Galois group $\mathcal{G} = \Gal(\mathcal{L}/
\mathcal{F})$ the set of all crossed homomorphisms $f:\mathcal{G}\to
\mathcal{L}^\times$ is, in fact, an Abelian group under the standard
operation of multiplying functions: for $f_1$ and $f_2$ crossed
homomorphisms from $\mathcal{G}$ to $\mathcal{L}^\times$, we define
$f_1\cdot f_2$ by $(f_1\cdot f_2)(g)=f_1(g)\cdot f_2(g)$ for $g$ in
$G$.  An important subgroup of the group of crossed homomorphisms is
the subgroup of \emph{coboundaries}.  A map $f: \mathcal{G}\to
\mathcal{L}^\times$ given by the formula $f(g) = {\lambda}/
{g(\lambda)}$ for some fixed $\lambda$ in $\mathcal{L}^\times$ is
called a \emph{coboundary}. Such a map is easily seen to be a
crossed homomorphism.  The first cohomology group $H^1 (\mathcal{G},
\mathcal{L}^\times)$ of $\mathcal{G}$ with coefficients in
$\mathcal{L}^ \times$ is defined to be the group of crossed
homomorphisms modulo the subgroup of coboundaries:
\begin{equation*}
    H^1 (\mathcal{G},\mathcal{L}^\times) = \frac{\{\mbox{ crossed
    homomorphisms }\}} {\{ \ \mbox{coboundaries}\ \}\ }.
\end{equation*}

In this language, Theorem~\ref{th:quadratic H90} translates to
$H^1(G, L^\times) = \{1\}$. Such an elegant outcome is what Speiser
sensed in the seedling: he proved that $H^1(\mathcal{G},
\mathcal{L}^\times) = \{1\}$ for all Galois extensions
$\mathcal{L}/\mathcal{F}$ with group $\mathcal{G}$.

\noindent \textbf{3. ``THE STALK IN TIME UNBENDS'': BIQUADRATIC
EXTENSIONS.} Replacing $G$ with the Klein $4$-group $\Z/2\Z\times
\Z/2\Z$, we shall observe the shoots of our seed ``bending double''
in this new context. It will appear, though, as if the new growth
points toward something unexpected, something akin to \emph{both}
the traditional Hilbert 90 result \emph{and} Speiser's cohomological
theorem on $H^1(G,L^\times)$.

Let $E$ be a Galois extension of $F$ with $G = \Gal (E/F) = \Z/2\Z
\times \Z/2\Z$. Since we assume that $F$ has characteristic other
than two, there exist $a_1$ and $a_2$ in $F^\times\setminus
(F^\times)^2$ such that $E = F (\sqrt{a_1}, \sqrt{a_2})$, and the
group $G$ is generated by two automorphisms $\sigma_1$ and
$\sigma_2$ such that
\begin{equation}\label{eq:deltaijcond}
    \frac{\sigma_i (\sqrt{a_j})}{\sqrt{a_j}}
    = (-1)^{\delta_{ij}},\tag{I}
\end{equation}
where $\delta_{ij}$ is Kronecker's delta function: $\delta_{ij}=0$
if $i=j$ and $\delta_{ij}=1$ if $i\neq j$. The lattice of subfields
of $E$ then takes the form
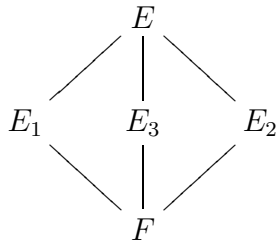
\begin{equation*}
    \xymatrix{& E \ar@{-}[dl] \ar@{-}[d]
    \ar@{-}[dr] \\ E_1 \ar@{-}[dr] & E_3
    \ar@{-}[d] & E_2 \ar@{-}[dl] \\ & F}
\end{equation*}
where $E_1 = F (\sqrt{a_1})$, $E_2 = F (\sqrt{a_2})$, and $E_3 =
F(\sqrt{a_1 a_2})$.  We say that an extension $E/F$ satisfying these
properties is a \emph{biquadratic extension}.

Since $E=E_1(\sqrt{a_2})$, $E=E_2(\sqrt{a_1})$, and $E =
E_3(\sqrt{a_1})$ are all quadratic extensions of fields of
characteristic different from two, the field extensions $E/E_1$,
$E/E_2$, and $E/E_3$ are also Galois extensions, with respective
Galois groups $\Gal(E/E_1)=\{\id, \sigma_1\}$, $\Gal(E/E_2) = \{\id,
\sigma_2\}$, and $\Gal(E/E_3)=\{\id, \sigma_1\sigma_2\}$.

Given quadratic Hilbert~90 for the two extensions $E/E_1$ and
$E/E_2$, we demonstrate that $H^1 (G,E^\times) = \{ 1 \}$ is
equivalent to another condition, itself very much in the spirit of
the original Hilbert 90 statement.  In fact, we show that this
condition is the ``difference'' between Speiser's result for
biquadratic extensions and quadratic Hilbert 90, in the following
sense: by adding the new result to quadratic Hilbert 90, we recover
Speiser's result in the biquadratic case.

Some corollaries of this condition have been rediscovered several
times, but the connection with Hilbert~90 and the natural proof that
we present seem new.  To establish this equivalence, we rephrase the
condition $H^1(G,E^\times) = \{1\}$ in the language of elements of
the multiplicative group $E^\times$.

Consider a crossed homomorphism $f : G \to E^\times$. Since $f(\id)
= 1$ and $f(\sigma_1 \cdot \sigma_2) = (\sigma_1 (f (\sigma_2)))
\cdot f (\sigma_1)$ we see that $f$ is determined by the values
$\alpha_i = f (\sigma_i)$. We further observe that
\begin{equation*}
    1 = f (\sigma_i^2) = (\sigma_i (\alpha_i)) \cdot \alpha_i =
    N_{E/E_i} (\alpha_i) \quad\quad (i = 1,2)
\end{equation*}
and that
\begin{equation*}
   (\sigma_1 (\alpha_2)) \cdot \alpha_1 = f (\sigma_1 \cdot \sigma_2)=
   f (\sigma_2 \cdot \sigma_1) = (\sigma_2 (\alpha_1)) \cdot \alpha_2.
\end{equation*}
Conversely, it is easy to check that for any given elements
$\alpha_1$ and $\alpha_2$ in $E^\times$ such that $N_{E/E_i}
(\alpha_i) = 1$ and  $(\sigma_1 (\alpha_2)) \cdot \alpha_1 =
(\sigma_2 (\alpha_1))\cdot \alpha_2$ there exists a unique crossed
homomorphism $f : G \to E^\times$ such that $f (\sigma_i) =
\alpha_i$.

Because $H^1(G,E^\times) = \{1\}$, we also know that for a given
crossed homomorphism $f$ there exists $\beta$ in $E^\times$ such that
$f(g) = {\beta}/{g(\beta)}$ for each $g$ in $G$.  In particular,
\begin{equation*}
    \alpha_i = f(\sigma_i) = \frac{\beta}{\sigma_i(\beta)}.
\end{equation*}
Therefore the cohomological identity $H^1 (G,E^\times) = \{1\}$ can
be reformulated as follows:
\begin{theorem}\label{th:biquadratic cohomology}
    Let $E=F(\sqrt{a_1},\sqrt{a_2})$ be a Galois extension of $F$
    with $G=\Gal(E/F)=\Z/2\Z \times \Z/2\Z$, and let $\sigma_1$ and
    $\sigma_2$ be generators of $G$ that satisfy
    relation~\eqref{eq:deltaijcond}.

    Then arbitrary elements $\alpha_1$ and $\alpha_2$ of $E^\times$
    satisfy conditions
    \begin{enumerate}
        \item $N_{E/E_1} (\alpha_1) =  1 = N_{E/E_2}(\alpha_2)$
    \end{enumerate}
    and
    \begin{enumerate}
        \setcounter{enumi}{1}
        \item\label{it:compcond} $\alpha_1\cdot \sigma_1(\alpha_2) =
        \alpha_2 \cdot \sigma_2 (\alpha_1)$
    \end{enumerate}
    if and only if
    \begin{enumerate}
        \setcounter{enumi}{2}
        \item there exists $\beta$ in $E^\times$ such that
        $\alpha_i = {\beta}/{\sigma_i(\beta)}$.
    \end{enumerate}
\end{theorem}
\noindent This result can be found in \cite[p.~756]{Co}.

Seen as results about kernels of norm functions,
Theorem~\ref{th:biquadratic cohomology} and Hilbert~90 are similar:
each shows that a kernel is suitably minimal. However, since
Theorem~\ref{th:biquadratic cohomology} describes the simultaneous
vanishing of norms from different quadratic extensions under the
additional compatibility condition~\eqref{it:compcond}, it is best
to think of Theorem~\ref{th:biquadratic cohomology} as a version of
Hilbert~90 with a compatibility condition.

We want to see just how far this new growth reaches above the
seed-case.  We find that, assuming quadratic Hilbert~90, the
conditions of Theorem~\ref{th:biquadratic cohomology} are equivalent
to a statement that the kernel of a particular operator is minimal.

In order to formulate this result, it is convenient to view $E^\times$
as a $\Z[G]$-module. The group ring $\Z[G]$ is defined to be the set
of formal sums of integer multiples of the elements of $G$,
\begin{equation*}
    \Z[G] = \{ c_0 \id + c_1 \sigma_1 + c_2 \sigma_2 + c_3 \sigma_1
    \sigma_2 : c_i \in \Z\},
\end{equation*}
together with addition and multiplication operations reminiscent of
the operations of polynomial addition and multiplication, as
follows. We first treat $\id$, $\sigma_1$, $\sigma_2$, and
$\sigma_1\sigma_2$ as indeterminates and add or multiply two
elements of $\Z[G]$ formally.  Then we replace products of these
``indeterminates'' with the corresponding product in the group $G$.
We derive the following formulas:
\begin{align*}
    (c_0 \id + &c_1 \sigma_1 + c_2 \sigma_2 + c_3 \sigma_1
    \sigma_2) + (d_0 \id + d_1 \sigma_1 + d_2 \sigma_2 + d_3 \sigma_1
    \sigma_2) := \\
    &(c_0+d_0)\id + (c_1+d_1)\sigma_1 + (c_2+d_2)\sigma_2 +
    (c_3+d_3)\sigma_1\sigma_2
\end{align*}
and
\begin{align*}
    (c_0 \id + &c_1 \sigma_1 + c_2 \sigma_2 + c_3 \sigma_1 \sigma_2)
    \cdot (d_0 \id + d_1 \sigma_1 + d_2 \sigma_2 + d_3 \sigma_1
    \sigma_2) := \\ &(c_0d_0+c_1d_1+c_2d_2+c_3d_3)\id +
    (c_0d_1+c_1d_0+c_2d_3+c_3d_2)\sigma_1 +\\
    &(c_0d_2+c_1d_3+c_2d_0+c_3d_1)\sigma_2 +
    (c_0d_3+c_1d_2+c_2d_1+c_3d_0)\sigma_1\sigma_2.
\end{align*}
A \emph{$\Z[G]$-module} $M$ is defined to be an Abelian group $M$
together with a given action of $\Z[G]$ on $M$.

In particular, the Abelian group $E^\times$ is a $\Z[G]$-module. We
specify the action of $\Z[G]$ on $E^\times$ by extending the action
of $\sigma_i$ on $E^\times$ in a natural way: for $\gamma$ in
$E^\times$ and $c_0,c_1,c_2,c_3$ in $\Z$ we set
\begin{equation*}
    (c_0\id + c_1 \sigma_1 + c_2 \sigma_2 + c_3 \sigma_1 \sigma_2)
    \cdot \gamma =
    \gamma^{c_0}\cdot \sigma_1(\gamma)^{c_1}\cdot
    \sigma_2(\gamma)^{c_2} \cdot \sigma_1\sigma_2(\gamma)^{c_3}.
\end{equation*}
Since all of the $\Z[G]$-modules we consider are, like $E^\times$,
multiplicative groups, we continue to write the module operation
multiplicatively. For example, $(\id+\sigma_1)\cdot m$ can be also
be written $m\cdot \sigma_1(m)$.  Similarly, we can write
$(\id-\sigma_1)\cdot m$ as $m/\sigma_1(m)$.  For ease of notation,
from now on we shall abbreviate $\id$ by $1$.  Thus $1$, considered
as an element of our group ring $\Z[G]$, will denote the identity
element.  From the context it will be clear when $1$ denotes an
element of $\Z[G]$. For instance, $(1+\sigma_i)\cdot m$ will denote
$(\id+\sigma_i)\cdot m$ for an element $m$ of a $\Z[G]$-module $M$.

Our approach to the new growth begins with a general theorem about
$\Z[G]$-modules for groups $G$ isomorphic to the Klein 4-group.  A
corollary that we deduce from the theorem will bring the new growth
into sharp relief.  In order to extend the property of quadratic
Hilbert 90 from $E^\times$ to general $M$, we say that
$(M,\sigma_1,\sigma_2)$ \emph{satisfies QH90} if for a given element
$m$ of $M$ the condition $(1+\sigma_i)\cdot m= 1$ implies that there
exists $n$ in $M$ such that $m=(1-\sigma_i)\cdot n$.

\begin{theorem}\label{th:newtheorem}
    Let $G=\Z/2\Z\times \Z/2\Z$, and let $\sigma_1$ and $\sigma_2$ be
    generators of $G$.  Let $M$ be a $\Z[G]$-module, written
    multiplicatively, and suppose that $(M,\sigma_1, \sigma_2)$
    satisfies QH90.  Then
    \begin{equation}\label{eq:star}
         \ker (1-\sigma_1) (1-\sigma_2) = \ker (1-\sigma_1) \cdot \ker
        (1-\sigma_2)\tag{II}
    \end{equation}
    if and only if for $m_1$ and $m_2$ in $M$ the two conditions
    \begin{enumerate}
      \setcounter{enumi}{0}
      \item\label{it:cond3} $(1+\sigma_1)\cdot m_1
      =1=(1+\sigma_2)\cdot m_2$
    \end{enumerate}
    and
    \begin{enumerate}
      \setcounter{enumi}{1}
      \item\label{it:cond4} $m_1\cdot\sigma_1(m_2)=
      m_2\cdot\sigma_2(m_1)$
    \end{enumerate}
    imply the truth of
    \begin{enumerate}
        \setcounter{enumi}{2}
        \item\label{it:cond5} there exists $n$ in $M$ such that
        $m_i=(1-\sigma_i)\cdot n$\quad $(i=1, 2)$.
    \end{enumerate}
\end{theorem}
\noindent In relation~\eqref{eq:star} the expression $\ker
(1-\sigma_1) \cdot \ker (1-\sigma_2)$ denotes the set of all
products $k_1\cdot k_2$ of elements $k_1$ in $\ker (1-\sigma_1)$ and
$k_2$ in $\ker (1-\sigma_2)$.

\begin{proof}
    Suppose first that
    \begin{equation*}
        \ker (1-\sigma_1)(1-\sigma_2) = \ker (1-\sigma_1)\cdot \ker
        (1-\sigma_2)
    \end{equation*}
    and that elements $m_1$ and $m_2$ of $M$ satisfy
    conditions~(\ref{it:cond3}) and (\ref{it:cond4}).  Since
    $(1+\sigma_i)\cdot m_i=1$, the QH90 property gives us elements
    $n_i$ in $M$ with $m_i=(1-\sigma_i)\cdot n_i$.
    Using~(\ref{it:cond4}) we obtain
    \begin{equation*}
        1 = \frac{(1-\sigma_2) \cdot m_1}{(1- \sigma_1)\cdot m_2} =
        (1-\sigma_2) (1-\sigma_1) \cdot (n_1/n_2).
    \end{equation*}
    Hence $n_1/n_2$ lies in $\ker (1-\sigma_1) (1-\sigma_2) = \ker
    (1-\sigma_1) \cdot \ker (1-\sigma_2)$. Recall that for each $i$
    the set $\ker (1-\sigma_i)$ is a group.  Therefore, if $e$ lies
    in $\ker (1-\sigma_i)$, then so does $1/e$.  Thus there exists
    $e_i$ in $\ker (1-\sigma_i)$ such that $n_1/n_2 = e_2/e_1$, and
    we set $n = n_1 \cdot e_1 = n_2 \cdot e_2$. We then see that
    \begin{equation*}
        (1-\sigma_i) \cdot n = m_i \quad\quad
        (i = 1,2),
    \end{equation*}
    establishing (\ref{it:cond5}).

    Now suppose that conditions (\ref{it:cond3}) and (\ref{it:cond4})
    imply condition (\ref{it:cond5}).  If $x = x_1 \cdot x_2$, where
    $x_i$ belongs to $\ker (1-\sigma_i)$, then
    \begin{align*}
        (1-\sigma_1) (1-\sigma_2) \cdot x &=
        \big((1-\sigma_2) (1-\sigma_1) \cdot x_1\big)\cdot \big(
        (1-\sigma_1) (1-\sigma_2) \cdot x_2\big) \\ &= 1,
    \end{align*}
    providing the inclusion $\ker(1-\sigma_1)\cdot \ker(1-\sigma_2)
    \subseteq \ker(1-\sigma_1)(1-\sigma_2)$.

    To verify the more interesting, opposite inclusion, consider an
    arbitrary $x$ in $\ker (1-\sigma_1) (1-\sigma_2)$. If $m_1 = (1-
    \sigma_1) \cdot x$ and $m_2 = 1$, then
    \begin{equation*}
        (1+\sigma_1) \cdot m_1 = \big((1+\sigma_1)(1-\sigma_1)\big)
        \cdot x = 1 = (1+\sigma_2) \cdot m_2.
    \end{equation*}
    Also
    \begin{equation*}
        \frac{(1-\sigma_2) \cdot m_1}{(1-\sigma_1) \cdot m_2} =
        (1-\sigma_2) (1-\sigma_1) \cdot x = 1.
    \end{equation*}
    Since we have verified (\ref{it:cond3}) and (\ref{it:cond4}) for
    $m_1$ and $m_2$, (\ref{it:cond5}) tells us that there exists $n$
    in $M$ such that $m_1 = (1-\sigma_1) \cdot x = (1-\sigma_1)
    \cdot n$ and $1 = m_2 = (1-\sigma_2) \cdot n$. Therefore $x / n$
    belongs to $\ker(1-\sigma_1)$ and $n$ to $\ker (1-\sigma_2)$.
    We conclude that $x = (x/n) \cdot n$ lies in $\ker (1-\sigma_1)
    \cdot \ker (1-\sigma_2)$.
\end{proof}

\begin{remark*}
    Since $\ker(1-\sigma_1) \cdot \ker(1-\sigma_2) \subseteq
    \ker(1-\sigma_1)(1-\sigma_2)$ for all $\Z[G]$-modules,
    relation~\eqref{eq:star} says that
    $\ker(1-\sigma_1)(1-\sigma_2)$ is as small as possible.
\end{remark*}

Returning to the $\Z[G]$-module $E^\times$, we observe that
$(E^\times,\sigma_1,\sigma_2)$ satisfies QH90. Indeed, for each $i$
and for $\gamma$ in $E^\times$ the condition $\gamma\cdot
\sigma_i(\gamma)=1$ implies that there exists $\delta$ in $E^\times$
such that $\gamma = \delta/\sigma_i(\delta)$.  Therefore
Theorems~\ref{th:biquadratic cohomology} and \ref{th:newtheorem}
applied to $M=E^\times$ immediately give the following corollary,
which makes the new growth visible:

\begin{corollary}\label{co:co1}
    Let $E=F(\sqrt{a_1},\sqrt{a_2})$ be a Galois extension of $F$
    with $G=\Gal(E/F)=\Z/2\Z \times \Z/2\Z$, and let $\sigma_1$ and
    $\sigma_2$ be generators of $G$ that satisfy
    relation~\eqref{eq:deltaijcond}.  Then
    \begin{equation*}
        \ker (1-\sigma_1) (1-\sigma_2) = \ker (1-\sigma_1) \cdot \ker
        (1- \sigma_2).
    \end{equation*}
\end{corollary}

\begin{remark*}
    By requiring elements $\sigma_1$ and $\sigma_2$ of $G$ to satisfy
    relation~\eqref{eq:deltaijcond}, we are viewing this result from
    only one particular direction. The equivalence in
    Theorem~\ref{th:biquadratic cohomology} implies the truth of every
    version of relation~\eqref{eq:star} obtained by substituting two
    arbitrary generators of $G$ in place of $\sigma_1$ and $\sigma_2$,
    while the conjunction of any version of relation~\eqref{eq:star}
    with quadratic Hilbert 90 gives the equivalence in
    Theorem~\ref{th:biquadratic cohomology}. Hence the particular
    ``difference'' we find in the corollary depends only on our vantage
    point.
\end{remark*}

\noindent\textbf{4. ``AND STARTS TO RAMIFY'': A FORM REVEALED.} Now
we establish relation~\eqref{eq:star} directly and record several of
its consequences. One corollary in particular has been rediscovered
several times under a certain attractive disguise as a statement on
binary quadratic forms.  Although several ingenious and beautiful
proofs of this statement have been obtained in the literature, the
proof we offer may be the most transparent.

\begin{theorem}\label{th:fullequiv}
    Let $E=F(\sqrt{a_1},\sqrt{a_2})$ be a Galois extension of $F$
    with $G=\Gal(E/F)=\Z/2\Z \times \Z/2\Z$, and let $\sigma_1$ and
    $\sigma_2$ be generators of $G$ that satisfy
    relation~\eqref{eq:deltaijcond}.  Then
    the following sets are identical:
    \begin{align*}
        K_1 &= \ker (1-\sigma_1) (1-\sigma_2),\\
        K_2 &= \ker (1-\sigma_1) \cdot \ker (1-\sigma_2),\\
        K_3 &= \langle E_1^\times,E_2^\times \rangle,\\
        K_4 &= \{ e \in E^\times : N_{E/E_3} (e) \in F^\times \},\\
        K_5 &= \{ e \in E^\times : N_{E/E_3} (e) \in
           N_{E_1/F}(E_1^\times) \cdot N_{E_2/F}(E_2^\times)\}.
    \end{align*}
\end{theorem}
\noindent In the definition of $K_3$, $\langle E_1^\times,
E_2^\times \rangle$ denotes the smallest subgroup of $E^\times$
containing $E_1^\times$ and $E_2^\times$, where $E_i=F(\sqrt{a_i})$.
In the definitions of $K_4$ and $K_5$, we recall that
$E_3=F(\sqrt{a_1a_2})$.  Finally, in the definition of $K_5$, the
product of norm groups denotes the set of all products of elements
of the first set with the second.

Taking $N_{E/E_3}$ of the sets $K_4$ and $K_5$ we obtain
\begin{equation*}
    N_{E/E_3} (E^\times) \cap F^\times \subseteq N_{E_1/F}
    (E_1^\times) \cdot N_{E_2/F} (E_2^\times).
\end{equation*}
The reverse inequality follows from the proof of
Theorem~\ref{th:fullequiv}, since for arbitrary $\gamma_i$ in
$E_i^\times$ and $\gamma = \gamma_1 \cdot \gamma_2$ we
have $N_{E/E_3} (\gamma) = N_{E_1/F}(\gamma_1) \cdot
N_{E_2/F}(\gamma_2)$. Thus we see that
\begin{equation}\label{eq:befcor}
    N_{E/E_3} (E^\times) \cap F^\times = N_{E_1/F} (E_1^\times)
    \cdot N_{E_2/F} (E_2^\times).\tag{III}
\end{equation}
This equality gives us, in the form of Corollary~\ref{co:co2}, the
equality of binary quadratic forms mentioned at the beginning of the
section. The proof of the corollary, which we present after the
proof of Theorem~\ref{th:fullequiv}, is a routine translation of
equation~\eqref{eq:befcor} into the language of binary quadratic
forms.
\begin{corollary}\label{co:co2}
    Suppose that $a$ and $b$ belong to $F^\times$ and $x$ and $y$ to
    $F (\sqrt{b})$. Then $x^2 - a y^2$ is in $F$ if and only if
    there exist $x_i$ and $y_i$ in $F$ $(i = 1,2)$ with the property
    that
    \begin{equation*}
        x^2 - a y^2 = (x_1^2 - a y_1^2) (x_2^2 - ab
        y_2^2).
    \end{equation*}
\end{corollary}
\noindent For a nice proof of and further references to this
corollary see \cite[Proposition~1.5]{LeSm} and the comments
preceding it.

\begin{proof}[Proof of Theorem~\ref{th:fullequiv}]
    $K_2=K_3$: First notice that $E = E_1(\sqrt{a_2}) =
    E_2(\sqrt{a_1})$.  Since the elements of $E$ fixed by the
    elements of $\Gal(E/E_i)$ are precisely those in $E_i$, we
    infer that $(1-\sigma_i)\cdot \ell =\ell/\sigma_i(\ell)=1$ if
    and only if $\ell$ lies in $E_i$. Since the domain of $\ker
    (1-\sigma_i)$ is $E^\times$, we obtain $\ker (1-\sigma_i) =
    E_i^\times$.

    \noindent $K_1\subseteq K_4$: Recall
    that
    \begin{equation*}
        \ker (1-\sigma_1) (1-\sigma_2) = \{ e \in E^\times :
         e \cdot \sigma_1 \sigma_2 (e) = \sigma_1 (e) \cdot
        \sigma_2 (e) \}.
    \end{equation*}
    Also $e \cdot \sigma_1 \sigma_2 (e) = N_{E/E_3} (e)$, a member
    of $E_3^\times$. Now assume that
    \begin{equation*}
        f = e \cdot \sigma_1 \sigma_2 (e) = \sigma_1 (e) \cdot
        \sigma_2 (e).
    \end{equation*}
    Then $\sigma_i (e\cdot \sigma_1 \sigma_2 (e)) = \sigma_1 (e)
    \cdot \sigma_2 (e) = f$. Because $\sigma_1$ and $\sigma_2$
    generate $G$, we see that $f$ is fixed by $G$.  Since the
    elements of $E$ fixed by $G$ lie in $F$, we deduce that $f$ lies
    in $F^\times$. Hence
    \begin{equation*}
      \ker (1-\sigma_1)(1-\sigma_2) \subseteq \{ e
      \in E^\times : N_{E/E_3} (e) \in F^\times \}.
    \end{equation*}

    \noindent $K_4\subseteq K_3$: Let $e$ in $E^\times$ have
    $f=N_{E/E_3}(e)$ in $F^\times$.  Because the set
    $\{1,\sqrt{a_1},\sqrt{a_2},\sqrt{a_1a_2}\}$ is basis for the
    vector space $E$ over the field $F$, we can write $e$ as a
    linear combination of the elements of this basis: $e=f_0 + f_1
    \sqrt{a_1} + f_2 \sqrt{a_2} + f_3 \sqrt{a_1 a_2}$. We calculate
    $N_{E/E_3}(e)$:
    \begin{align*}
        f &= N_{E/E_3} (e) = e \cdot \sigma_1\sigma_2(e) \\ &=
        e\cdot (f_0 - f_1 \sqrt {a_1} - f_2 \sqrt{a_2} + f_3
        \sqrt{a_1 a_2}).
    \end{align*}
    Multiplying out the product and determining the coefficient of
    the basis element $\sqrt{a_1a_2}$, we learn that $f_0 f_3 - f_1
    f_2 = 0$.

    Assume first that both $f_2$ and $f_3$ are nonzero.  Thus $f_0 =
    f_2 t$ and $f_1 = f_3 t$ for some $t$ in $F^\times$.
    Substituting these values into our expression for $e$ gives
    \begin{equation*}
        e = (f_2 + f_3 \sqrt{ a_1}) \cdot (t + \sqrt{a_2})  \in
        \langle E_1^\times, E_2^\times \rangle,
    \end{equation*}
    the subgroup of $E^\times$ generated by $E_1^\times$ and
    $E_2^\times$. If $f_2 = 0$, then $f_0 = 0$ or $f_3 = 0$. In the
    first case
    \begin{equation*}
        e = \sqrt{a_1} \cdot \left(f_1 + f_3 \sqrt{a_2}\right) \in
        \langle E_1^\times, E_2^\times \rangle,
    \end{equation*}
    and in the second case
    \begin{equation*}
        e = f_0 + f_1 \sqrt{a_1} \in E_1^\times \subseteq \langle
        E_1^\times, E_2^\times \rangle.
    \end{equation*}
    The case $f_3 = 0$ is handled in the same way. Thus we see that
    \begin{equation*}
        \{ e \in E^\times : N_{E/E_3} (e) \in F^\times \} \subseteq
        \langle E_1^\times, E_2^\times \rangle.
    \end{equation*}

    \noindent $K_3\subseteq K_1$: Since the elements of $E$ fixed by
    $\Gal(E/E_i)$ lie in $E_i$, we find that $E_i^\times \subseteq
    \ker (1-\sigma_i)$ $(i=1,2)$.  We therefore obtain
    \begin{equation*}
      E_i^\times \subseteq \ker (1-\sigma_1)(1-\sigma_2),
    \end{equation*}
    from which the stated inclusion follows.

    We have established that the sets $K_1$, $K_2$, $K_3$, and $K_4$
    coincide, so it remains to show that set $K_5$ is identical to
    the others.  The inclusion $K_5\subseteq K_4$ is easy to see:
    \begin{align*}
        W &:= \{ e \in E^\times :  N_{E/E_3} (e) \in
        N_{E_1/F} (E_1^\times) \cdot N_{E_2/ F} (E_2^\times) \} \\
        &\subseteq\ \{ e \in E^ \times :  N_{E/E_3} (e) \in
        F^\times \}.
    \end{align*}
    We argue that $K_3\subseteq K_5$, as follows.  Let
    $\gamma = \gamma_1 \cdot \gamma_2$, where $\gamma_i$
    is in $K_i^\times$. Note that $N_{E/E_3}(\gamma_i) =
    \gamma_i \cdot \sigma_1 \sigma_2 (\gamma_i) = N_{E_i/F}
    (\gamma_i)$ because $\sigma_1 \sigma_2$ acts on both $E_1$ and
    $E_2$ nontrivially. It follows that
    \begin{equation*}
        N_{E/E_3}(\gamma) \in N_{E_1/F} (E_1^ \times) \cdot N_{E_2/F}
        (E_2^\times).
    \end{equation*}
    Therefore $\gamma$ belongs to $W$. Hence the sets $K_1$, $K_2$,
    $K_3$, $K_4$, and $K_5$ are identical.
\end{proof}

\begin{proof}[Proof of Corollary~\ref{co:co2}]
    Recall that the nonzero values of a binary quadratic form $x^2 -
    d y^2$ for $d$ in $F$ form a group under multiplication. Hence
    if $b$ is a square in $F$ our statement follows. Therefore
    assume that $b$ lies in $F^\times \setminus (F^{\times})^2$.

    Next observe that the set of values assumed by a quadratic form
    $x^2 - c^2 y^2$ $(c \ne 0)$ as $x$ and $y$ range over $F$ is
    $F$ itself:
    \begin{equation*}
        f = \left(\left(\frac{f+1}{2}\right)^2 -
        c^2 \left(\frac{f-1}{2c}\right)^2\right) \quad\quad
        (f\in F).
    \end{equation*}
    Therefore if either $a$ or $ab$ is a square in $F$, the statement is
    true as well. (Notice that if $ab$ is a square in $F$, then $a$
    is a square in $F (\sqrt{b})$ and $x^2 - a y^2$ assumes
    all elements of $F(\sqrt{b})$.)

    Finally, assume that none of the elements $a$, $b$, or $ab$ is a
    square in $F$. Let $a_1 = a$ and $b = a_2/a_1$. Since $E = E_3
    (\sqrt{a_1}) = E_3 (\sqrt{a})$, we find that
    \begin{equation*}
        N_{E/E_3} (E^\times) \cap F^\times = \{ x^2 - a y^2 :
        x,y \in E_3, \ x^2 - a y^2 \in F^\times \}.
    \end{equation*}
    Also note that $E_3 = F (\sqrt{a_1 a_2}) = F (a_1 \cdot \sqrt{
    a_2/a_1}) = F (\sqrt{b})$. From $N_{E_i/F} (E_i^ \times) = \{
    x_i^2 - a_i y_i^2 : x_i, y_i \in F, \ x_i^2 - a_i y_i^2 \ne 0
    \}$ the result follows.
\end{proof}

\noindent\textbf{5. ``BE VAGUELY VAST, AND CLIMB'': FURTHER
DIRECTIONS.} We have observed the first shoots of the generalization
from quadratic to biquadratic extensions.  To understand what truly
lies cached within the seed-case, we must watch closely as $E/F$ is
replaced with bicyclic extensions, tricyclic extensions, and even
general Abelian extensions. In particular, it is natural to ask what
the ``difference'' is between bicyclic Hilbert 90 and biquadratic
Hilbert 90, and, more generally, between Hilbert 90 for an Abelian
extension and Hilbert 90 for a bicyclic extension. Perhaps a
generalization of Theorem~\ref{th:fullequiv} similarly connects a
notion of cohomological vanishing with minimal kernels of related
operators.

Like seedlings ourselves, however, we have been ``forced to make
choice of ends,'' and have accepted---not unlike this article
itself---``the doom of taking shape.''  But what you, the reader,
may yet heave aloft---``that remains,'' as the poet says, ``to be
seen.''

\medskip

\noindent\Small{\textbf{ACKNOWLEDGMENTS.}\ \  We thank Noam Elkies
for bringing to our attention the reference \cite{O}, which in turn
had been brought to his attention by Franz Lemmermeyer.  We also
thank both referees for thoughtful comments and valuable suggestions
that have improved our exposition.  Roman Dwilewicz acknowledges
support from the University of Missouri Research Board and the
Polish Committee for Scientific Research grant KBN2P03A04415. J\'an
Min\'a\v{c} is grateful for support from Natural Sciences and
Engineering Research Council of Canada grant R0370A01 and a
Distinguished Professorship for 2004--2005 at the University of
Western Ontario. Andrew Schultz thanks Ravi Vakil for his
encouragement and direction in this and all other projects.  John
Swallow thanks the Institut de Math\'ematiques de Bordeaux for its
hospitality during part of the writing of this article.}


\parskip=6pt plus 7pt

\Small
\noindent{\textbf{ROMAN DWILEWICZ} is indigenous to Poland.  He
earned his master's, doctoral, and habilitation degrees at the
University of Warsaw, and then for many years he taught at the same
institution. He has been a member of the Polish Academy of Sciences,
as well as scientific director of its Institute of Mathematics and
the Banach Center.  Once transplanted halfway around the world, he
met J\'an Min\'a\v{c} when J\'an, absorbed in thought, quite
literally crashed into him.  \\
\textit{Department of Mathematics and Statistics, University of
Missouri-Rolla, Rolla, MO 65409}\\
\textit{romand{\@@}umr.edu}}

\noindent{\textbf{J\'AN MIN\'A\v{C}} spent years blissfully unaware
of the fact that the Fields Medal clock had already begun ticking
from the moment of his birth. He happily wasted his time
daydreaming, playing soccer, reading all kinds of literature, and
toying with prime numbers. Although now past forty, he is as excited
as ever to keep playing with mathematics, blessed with great
collaborators and the support of his wonderful wife Leslie. Summing
these up, he reaches but one result, his own sort of Fields Medal.
\\
\textit{Department of Mathematics, Middlesex College,
University of Western Ontario, London, ON N6A 5B7, Canada}\\
\textit{minac{\@@}uwo.ca}}

\noindent{\textbf{ANDREW SCHULTZ} arrived at Davidson College
intending to study economics.  What was to be a brief fling with
mathematics in his freshman year, however, quickly blossomed into
love, and he chose to abandon all practicality in its pursuit. He
spent the summers of 2001 and 2002 in the REU run by Arie
Bialostocki at the University of Idaho, first as a participant and
then as assistant director. He is currently a graduate student at
Stanford, studying algebraic geometry with Ravi Vakil.
\\
\textit{Department of Mathematics, Stanford University, Stanford, CA
94305-2125}\\
\textit{aschultz{\@@}stanford.edu}}

\noindent{\textbf{JOHN SWALLOW} is Kimbrough Associate Professor at
Davidson College. In moments of pure serendipity he met, in the same
semester, a sophomore named Andrew Schultz and a fellow MSRI member
named J\'an Min\'a\v{c}. He is the author of \textit{Exploratory
Galois Theory} (Cambridge, 2004) and is the 2006--2007 section
lecturer for the MAA's Southeastern Section. He has been elected a
trustee of his alma mater, Sewanee: The University of the South.
\\
\textit{Department of Mathematics, Davidson College, Box 7046,
Davidson, NC 28035}\\
\textit{joswallow{\@@}davidson.edu} }

\begin{thebibliography}{14}
\bibitem[1]{Co} I.~G.~Connell, Elementary generalizations of
Hilbert's Theorem~90, \textit{Canad. Math.~Bull.} \textbf{8}
(1965) 749--757.

\bibitem[2]{El} N.~D.~Elkies, Pythagorean triples and
Hilbert's Theorem~90, this {\sc Monthly} \textbf{110} (2003)
678.

\bibitem[3]{H} D.~Hilbert, \emph{The Theory of Algebraic Number
Fields} (trans.~I.~Adamson), Springer-Verlag, Berlin, 1998.

\bibitem[4]{Kum} E.~Kummer, \"Uber eine besondere Art aus
complexen Einheiten gebildeter Ausdr\"ucke, \emph{J.~Reine
Angew.~Math.} {\bf 50} (1855) 212--232.

\bibitem[5]{LeSm} D.~B.~Leep and T.~L.~Smith, Multiquadratic
extensions, rigid fields and Pythagorean fields, \emph{Bull. London
Math. Soc.} {\bf 34} (2002) 140--148.

\bibitem[6]{LeSc} F.~Lemmermeyer and N.~Schappacher, Introduction,
in
D.~Hilbert, \emph{The Theory of Algebraic Number Fields}
(trans.~I.~Adamson), Springer-Verlag, Berlin, 1998, pp.~xxiii--xxxvi.

\bibitem[7]{N} E.~Noether, Der Hauptgeschlechtssatz f\"ur
relativ-galoissche Zahlk\"orper, \emph{Math. Ann.} {\bf 108} (1933)
411--419.

\bibitem[8]{O} T.~Ono, \emph{Variations on a Theme of Euler.
Quadratic Forms, Elliptic Curves, and Hopf Maps}, Plenum,
New York, 1994.

\bibitem[9]{S} A.~Speiser, Zahlentheoretische S\"atze aus der
Gruppentheorie, \emph{Math.~Zeit.} {\bf 5} (1919) 1--6.

\bibitem[10]{T} O.~Taussky, Sums of squares, this {\sc Monthly}
{\bf 77} (1970) 805--830.

\bibitem[11]{Vo1} V.~Voevodsky, Motivic cohomology with
$\Z/2$-coefficients, \textit{Publ.~Inst.~Hautes \'Etudes Sci.}
{\bf 98} (2003) 59--104.

\bibitem[12]{Vo2} ---------, On motivic cohomology with $\Z/l$
coefficients, $K$-Theory Preprint Archive 639, 2003; available at
www.math.uiuc.edu/K-theory/0639/.

\bibitem[13]{We} E.~Weiss, \emph{Cohomology of Groups},
Academic Press, New York, 1969.

\bibitem[14]{Wi} R.~Wilbur, ``Seed Leaves,'' \emph{New and Collected
Poems}, Harcourt Brace Jo\-va\-no\-vich, San Diego, 1988, pp.~129--130.
\end{thebibliography}
\end{document}